\documentclass[12pt]{article}

\usepackage{amsthm}
\usepackage{latexsym}
\usepackage{amssymb,amsmath}
\usepackage{epsfig}
\usepackage{amsfonts}
\usepackage{graphicx}
\usepackage{color}
\usepackage{algorithm}
\usepackage{algorithmic}
\usepackage{caption}
\usepackage{url}
\usepackage{cite}

\date{}

\newcommand{\gp}{{\rm gp}}
\newcommand{\Fitness}{{\rm Fitness}}
\newcommand{\Cay}{{\rm Cay}}

\DeclareUnicodeCharacter{200B}{}
\DeclareUnicodeCharacter{041D}{}
\DeclareUnicodeCharacter{0430}{}
\DeclareUnicodeCharacter{0443}{}
\DeclareUnicodeCharacter{0447}{}
\title{Three algorithmic approaches to the general position problem}

\author{Zahra Hamed-Labbafian$^{a,}$\thanks{Email: \texttt{hamedlabbafianzahra@gmail.com}} 
\and  
Narjes Sabeghi$^{b,}$\thanks{Corresponding author, Email: \texttt{n.sabeghi@velayat.ac.ir}} 
\and 
Mostafa Tavakoli$^{a,}$\thanks{Email: \texttt{m$\_$tavakoli@um.ac.ir}} 
\and 
Sandi Klav\v{z}ar$^{c,d,e,}$\thanks{Email: \texttt{sandi.klavzar@fmf.uni-lj.si}}
}

\begin{document}

\maketitle

\begin{center}
$^a$ Department of Applied Mathematics,
Faculty of Mathematical Sciences \\
Ferdowsi University of Mashhad,
Mashhad, Iran\\
\medskip

$^b$ Department of Mathematics,
Faculty of Basic Sciences\\
Velayat University, Iranshar, Iran\\
\medskip

$^c$ Faculty of Mathematics and Physics, University of Ljubljana, Slovenia\\
\medskip

$^d$ Institute of Mathematics, Physics and Mechanics, Ljubljana, Slovenia\\
\medskip

$^e$ Faculty of Natural Sciences and Mathematics, University of Maribor, Slovenia
\end{center}

\date{}
\maketitle

\begin{abstract}
If $G$ is a graph, then $X\subseteq V(G)$ is a general position set if for every two vertices $v,u\in X$ and every shortest $(u,v)$-path $P$, it holds that no inner vertex of $P$ lies in $X$. In this note we propose three algorithms to compute a largest general position set in $G$: an integer linear programming algorithm, a genetic algorithm, and a simulated annealing algorithm. These approaches are supported by examples from different areas of graph theory.
\end{abstract}

\\

\noindent
{\bf AMS Subj.\ Class.\ (2020)}:  05C12, 05C69, 05C76, 05C85

\baselineskip16pt

\section{Introduction}

The general position problem in graph theory was independently introduced (at least) three times, which clearly demonstrates that the concept is of wider interest. The first time it was investigated in the context of hypercubes~\cite{korner-1995}, and later extended to the general case~\cite{chandran-2016, manuel-2018a}. The latter article has led to extensive research of the problem, and a recent review article~\cite{chandran-2025} lists 115 references. It should be emphasized that the study of general position sets has extended to various types, such as Steiner position sets~\cite{klavzar-2021}, edge general position sets~\cite{manuel-2022}, mutual-visibility sets \cite{distefano-2022}, mobile position sets~\cite{klavzar-2023}, monophonic position sets~\cite{thomas-2024}, vertex position sets~\cite{thankachy-2024}, and lower general position sets~\cite{distefano-2025}. Among the recent achievements in the field, we would like to highlight~
\cite{araujo-2025, irsic-2024, thomas-2024a, tian-2025, tian-2024, tuite-2025}.

Let's now define the problem. If $G = (V(G), E(G))$ is a simple, connected graph, then the set $X\subseteq V(G)$ is a {\em general position set} for the graph $G$ if for every two vertices $v,u\in X$ and every shortest $(u,v)$-path $P$, it holds that no inner vertex of $P$ lies in $X$. We are interested in the size of a largest general position set of $G$ which is called the \textit{general position number} of $G$ and denoted by $\gp(G)$.

The {\sc General Position Set} problem for a given graph $G$ and a positive integer $k\leq |V(G)|$ asks whether there exists a general position set $S$ for $G$ such that $|S|\geq k$. Already in the seminal paper~\cite{manuel-2018a} it was shown that the {\sc General Position Set} problem is NP-complete for arbitrary graphs. The closely related {\sc Lower General Position} problem was proved to be NP-complete in~\cite{distefano-2025}, and the {\sc Monophonic Position Set} problem was demonstrated to be NP-hard in~\cite{thomas-2024}. Interestingly, it is not clear whether the {\sc Monophonic Position Set} is NP-complete for general graphs.

Since determining the general position number of a graph is difficult in general, it is reasonable to approach it in different algorithmic ways. Kor\v{z}e and Vesel~\cite{korze-2023} described a reduction from the problem of finding a general position set in a graph to
the satisfiability problem and applied the approach to the case of general position sets of hypercubes. However, we are aware of no systematic algorithmic approach to the {\sc General Position Problem} so far. Therefore, in this note we propose three different algorithms for approaching the problem. These approaches are integer linear programming (ILP), genetic algorithm (GA), and simulated annealing (SA), and are described in the next section. In Section~\ref{sec:experiments} these approaches are tested on several large graphs from different areas of graph theory. 

\section{Three approaches}
\label{sec:three-approaches}

In this section, we present in respective subsections the three approaches to compute the general position number as announced in the introduction.  

\subsection{Integer linear programming model}

Let $G$ be a graph with $V(G) = \{v_1,\ldots, v_{n(G)}\}$, where $n(G)$ denotes the order of $G$. To find a general position set of $G$ with the maximum cardinality, for a given set $X\subseteq V(G)$ define binary variables $x_j$, $j\in [n(G)]$, which constitute a characteristic vector of $X$, that is, 
\begin{equation*}
x_j = \left\lbrace
\begin{array}{ll}
1; &  v_j \in X\,, \\
0; & \text{otherwise}\,.
\end{array} \right.
\end{equation*}
Then the objective function is to maximize $\sum_{j = 1}^{n} x_j$. For each two vertices $v_i$ and $v_j$, let $I_G(v_i,v_j)$ be the set of vertices different from $v_i$ and $v_j$ which lie on at least one shortest $(v_i,v_j)$-path. Therefore, if $v_i, v_j\in X$, that is, if $x_i + x_j = 2$, then we must have $\sum_{x_\ell \in I_G(v_i , v_j)} x_\ell = 0$, 
which can be considered as the inequality:
\[
\sum_{x_\ell \in I(v_i , v_j)}x_\ell \le 0\,.
\]
Using the integer modeling techniques, we model this constraint as:
 \begin{align*}
  \sum_{v_\ell \in I(v_i , v_j)}x_\ell + M(x_i + x_j-1) \leq M, \quad i,j \in [n(G)], 
  \end{align*}
where $M$ is an upper bound for the constraint $\sum_{v_\ell \in I(v_i , v_j)}x_\ell \leqslant 0$. In this case we can select $M = n(G)$. ILP for our problem is thus:
 \begin{align*}
  & \max \sum_{j = 1}^{n} x_j \\
   & \ {\rm s.t.} \quad \sum_{v_\ell \in I(v_i , v_j)}x_\ell + n(G)(x_i + x_j) \leq 2n(G), \quad i,j \in [n(G)] \\
 & \qquad \ \ x_i \in \{0,1\}, \quad i\in [n(G)]
  \end{align*}

\subsection{Genetic Algorithm}

In this subsection, we present a Genetic Algorithm (GA for short) to find largest general position sets. GA is a population-based method and is well-suited for the general position problem because it efficiently explores large search spaces while maintaining feasible solutions through its crossover and mutation operators.

Before implementing the GA, it is necessary to define the stopping criterion, the fitness function, and the genetic operators specific to this problem. Each instance of solution $S$ for a graph $G$ is represented as a binary (characteristic) vector of size $n(G)$, that is, if the vertex $i$ is included in the solution, then $S(i)=1$, otherwise, $S(i)=0$. The stopping criterion is defined as the maximum number of iterations. The objective of the algorithm is to obtain a solution that satisfies the general position constraints while maximizing the number of selected vertices. The fitness function for a solution $S$ is defined as: 
$$
\Fitness(S) = \sum S - M \cdot f,
$$
where $\sum S$ represents the total number of selected vertices, and $f$ denotes the number of violations of the general position condition. Specifically, for every pair of vertices $u$ and $v$ with $S(u) = 1$ and $S(v) = 1$, if there exists a vertex $w$ on a shortest $(u,v)$-path with $S(w) = 1$, then $f$ is increased. The penalty term $M \cdot f$ ensures that solutions violating the general position condition are penalized proportionally to the number of violations, where $M$ is a large constant balancing the trade-off between maximizing the number of selected vertices and maintaining feasibility.
The initial population is generated randomly, ensuring diversity in the search space. 

For parent selection, we use a random selection method. The offspring generation process involves combining two parent solutions: the new solution instance is a binary vector, where a vertex is included if it appears in at least one of the parents. 

The mutation operator selects two vertices randomly and swaps their values. 

The GA presented in Algorithm~\ref{alg1} outlines the general structure of the approach for solving the general position problem, while the details of the various operators have been explained above. 


\begin{algorithm}[ht!]
\caption{GA for the general position problem}\label{alg1}
\begin{algorithmic}[1]
\REQUIRE Graph, $G$, maximum number of iterations $Maxit$, and initial population of size $n_p$.
\ENSURE a general position set and $\gp(G)$.
\STATE Create a feasible initial populations with size $n_p$, and evaluate them using fitness function. 
\STATE $count=1$
\WHILE {$count\leq Maxit$}
\STATE Select parents randomly.
\STATE Generate and evaluate offspring using the crossover operator.
\STATE Generate and evaluate mutated populations using the mutation operator.
\STATE Merge the initial population, offspring, and mutated populations, and sort them based on the fitness function (new population).
\STATE Truncate new population and generate a new population with size $n_p$.
\STATE  $count=count+1$
\ENDWHILE
\STATE Return the individual with the best fitness value as the best solution.
\end{algorithmic}
\end{algorithm}

\subsection{Simulated Annealing}

Simulated Annealing (SA for short) is a single-solution point-based metaheuristic algorithm inspired by the annealing process in metallurgy. Before implementing SA for the general position problem, it is essential to define the stopping criterion, the neighborhood generation mechanism, and the fitness function. The algorithm begins with an initial solution $S_0$, which is a vector of zeros and ones representing the inclusion or exclusion of vertices in the solution set. The quality of a solution is measured by the fitness function, $\Fitness(S)$, which is defined in the same way as for the GA. 

The acceptance probability of a new solution is determined by the Metropolis criterion, which depends on the current temperature and the change in the objective value. The temperature gradually decreases according to a predefined cooling schedule, and the algorithm terminates once the maximum number of iterations or a minimum temperature is reached. 
The SA is formally presented in Algorithm~\ref{alg3}. 
\begin{algorithm}[ht!]
\caption{Simulate Annealing (SA) algorithm for the general position problem}\label{alg3}
\begin{algorithmic}[1]
\REQUIRE Graph $G$, maximum number of iterations $Maxit$, initial Temperature $T_0$, and the rate of temperature reduction $\rho$.
\ENSURE A general position set of $G$ and $\gp(G)$.
\STATE Create a feasible initial solution $S$.
\STATE $best\_solution\leftarrow S$
\STATE $best\_fitness \leftarrow S_{Fit}$ ($S_{Fit}=Fitness(S)$)
\STATE $count\leftarrow 1$
\STATE $T\leftarrow T_0$
\WHILE {$count\leq Maxit$}
\STATE Generate and evaluate Neigbors of $S$.
\STATE Let $best\_neighbor$ be the neighbor with the best fitness value ($best\_neighbor\_Fit$) 
\STATE $S\leftarrow best\_neighbor$, $S_{Fit}\leftarrow best\_neighbor\_Fit$
\IF {$S_{Fit}>best\_fitness$}
\STATE $best\_solution\leftarrow S$
\STATE $best\_fitness \leftarrow S_{Fit}(Fitness(S))$
\ELSE
\STATE Generate a random number $r$.
\IF {$exp(-S_{Fit})/T>r$}
\STATE $best\_solution\leftarrow S$
\STATE $best\_fitness \leftarrow S_{Fit}(Fitness(S))$
\ENDIF
\ENDIF
\IF {$k \geq cooling\_time$} 
\STATE $T\leftarrow T*\rho$
\STATE $k\leftarrow 0$
\ENDIF
\STATE $k\leftarrow k+1$
\STATE  $count\leftarrow count+1$
\ENDWHILE
\STATE Return the $best\_solution$ with the best fitness value as the best solution.
\end{algorithmic}
\end{algorithm}

\section{Experimental results}
\label{sec:experiments}

We have performed our experiments on three different types of graph. The first examples come from mathematics chemistry, the second examples deal with hypercubes, while the last examples come from algebraic graph theory.  

\begin{figure}[ht!]
\centerline{
 \includegraphics[scale=.5]{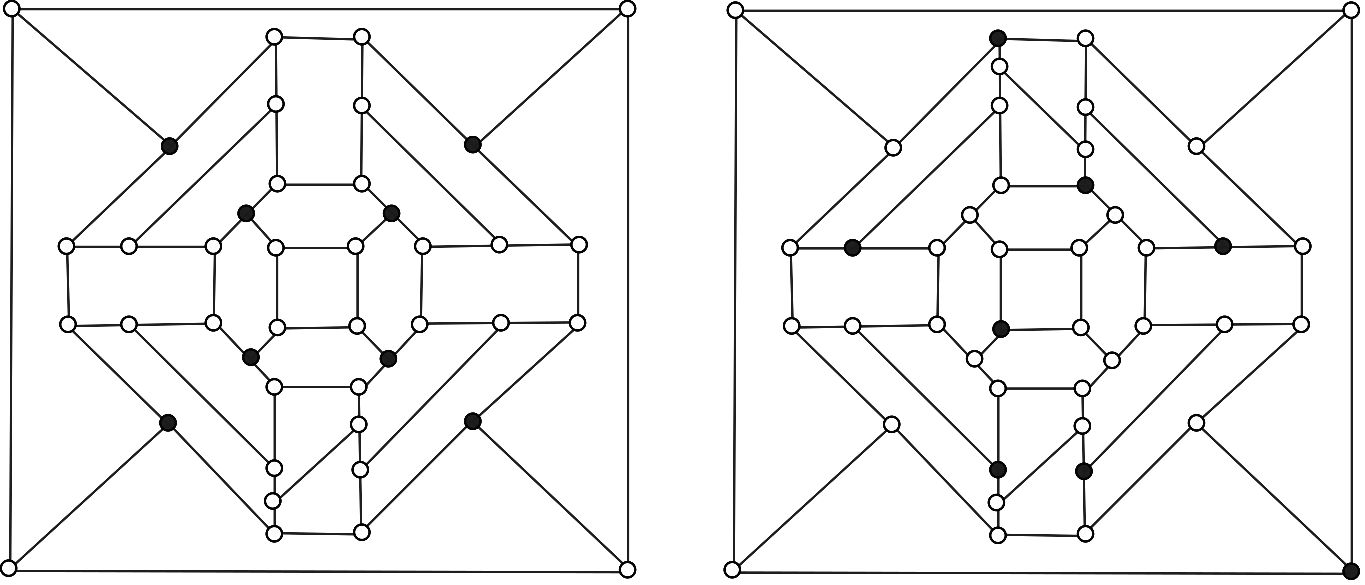}}
 \caption{Maximum general position sets for $C_{42}$ and $C_{44}$}\label{Fig2}
   \end{figure}

\begin{figure}[ht!]
\centerline{
 \includegraphics[scale=.5]{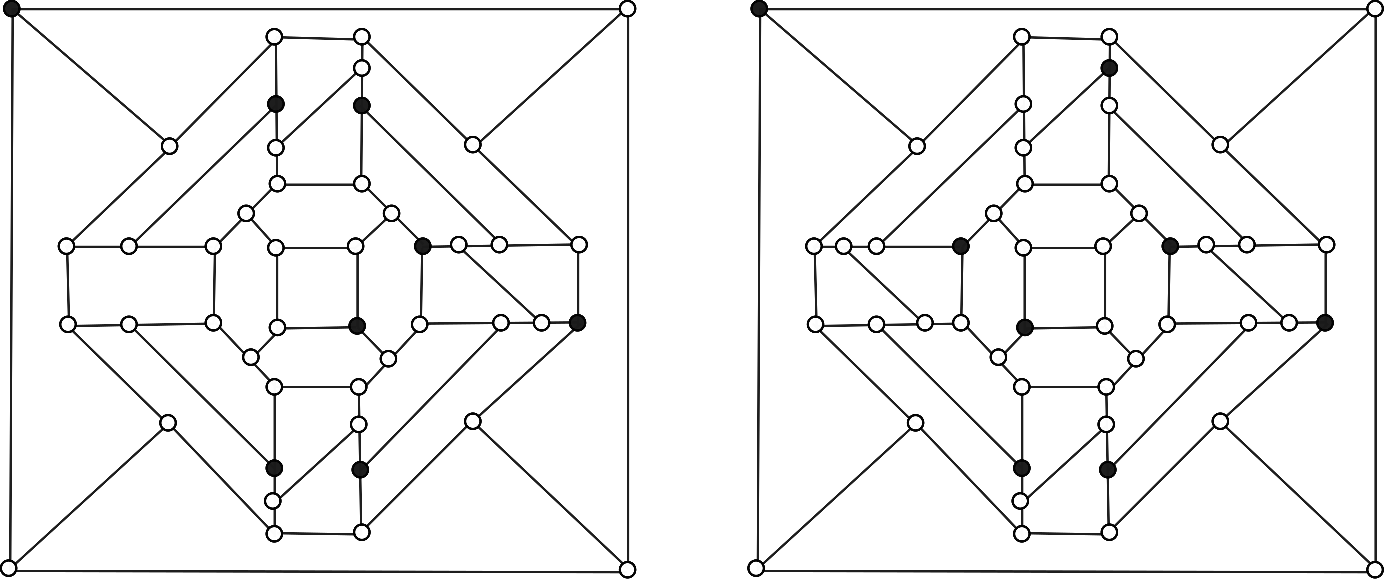}}
 \caption{Maximum general position sets for $C_{46}$ and $C_{48}$}\label{Fig3}
   \end{figure}

From the area of chemical graph theory we have considered the so-called middle-fullerene graphs $C_{42}$, $C_{44}$, $C_{46}$, and $C_{48}$ from~\cite{melker-2017}, see Figs.~\ref{Fig2} and~\ref{Fig3}. 

The maximum general position sets of $C_{42}$ and $C_{44}$ as indicated in Fig.~\ref{Fig2} with solid dots, were obtained by solving the corresponding ILP model. In the case of the middle-fullerene graphs $C_{46}$ and $C_{48}$, their examination was carried out using all three algorithms, the results are reported in the top two lines of Table~\ref{Tab1}, the concrete largest general position sets obtained can be seen in Fig.~\ref{Fig3}.

Let $Q_n$ denote the $n$-dimensional hypercube. It is known all the way from~\cite{korner-1995} that it is intrinsically difficult to determine $\gp(Q_n)$ in general. The small values $\gp(Q_1) = \gp(Q_2) = 2$, and $\gp(Q_3) = 4$ can be determined by hand, while the values $\gp(Q_4) = 5$, $\gp(Q_5) = 6$, $\gp(Q_6) = 8$, and $\gp(Q_7) = 9$ were reported in~\cite[Table~1]{korze-2023}. In this note we have verified these values by our three algorithms, the computational results are presented in the middle of Table~\ref{Tab1}. 

For the last set of examples we have considered three examples from algebraic graph theory, Cayley graphs
$\Cay(\mathbb{Z}_9,\{1,3,6,8\})$,  $\Cay(\mathbb{Z}_9,\{1,2,3,6,7,8\})$, 
and
$\Cay(\mathbb{Z}_{20},\{1,3,17,19\})$. The computational results are summarized in the last three rows of Table~\ref{Tab1}. 

\begin{table}[ht!]
\begin{center}
\begin{tabular}{|c|c|c|c|c|c|}
\hline
$G$ & $n(G)$  &\begin{tabular}[x]{@{}c@{}}${\gp}_{\rm ILP}$ \\ ${\gp}_{\rm GA}$ \\ ${\gp}_{\rm SA}$ \end{tabular}&ILP time&GA time&SA time\\ \hline
$C_{46}$&46 &8&73.10&\begin{tabular}[x]{@{}c@{}}29.72\\ $n_p=50 $\\$Maxit=5000$\end{tabular}&\begin{tabular}[x]{@{}c@{}}22.32\\ $T_0=10$\\$Maxit=500$\end{tabular} \\ \hline

$C_{48}$&48 &8&71.46&\begin{tabular}[x]{@{}c@{}}38.34\\ $n_p=50$\\$Maxit=5000$\end{tabular}&\begin{tabular}[x]{@{}c@{}}11.95\\ $T_0=10$\\$Maxit=500$\end{tabular}\\ \hline

$Q_{3}$&8 &4&0.03&\begin{tabular}[x]{@{}c@{}}0.05\\ $n_p=10$\\$Maxit=100$\end{tabular}&\begin{tabular}[x]{@{}c@{}}0.010\\ $T_0=10$\\$Maxit=10$\end{tabular}\\ \hline

$Q_{4}$&16 &5&0.38&\begin{tabular}[x]{@{}c@{}}0.22\\ $n_p=20$\\$Maxit=200$\end{tabular}&\begin{tabular}[x]{@{}c@{}}0.075\\ $T_0=10$\\$Maxit=10$\end{tabular}\\ \hline

$Q_{5}$&32 &6&11.52&\begin{tabular}[x]{@{}c@{}}1.84\\ $n_p=20$\\$Maxit=400$\end{tabular}&\begin{tabular}[x]{@{}c@{}}1.93\\ $T_0=10$\\$Maxit=50$\end{tabular}\\ \hline

$Q_{6}$&64 &8&359.01&\begin{tabular}[x]{@{}c@{}}82.11\\$n_p=50$\\$Maxit=4500$\end{tabular}&\begin{tabular}[x]{@{}c@{}}144.985\\ $T_0=10$\\$Maxit=500$\end{tabular}\\ \hline

$Q_{7}$&128 &9&3617.39&\begin{tabular}[x]{@{}c@{}}2060.06\\ $n_p=50$\\$Maxit=8000$\end{tabular} &\begin{tabular}[x]{@{}c@{}}2574.79 \\ $T_0= 10$\\$Maxit=100 $\end{tabular}\\ \hline

$\Cay(\mathbb{Z}_9,\{1,3,6,8\})$ &9&4&1.22&\begin{tabular}[x]{@{}c@{}}0.06\\
$n_p=10$\\$Maxit=50$\end{tabular}&\begin{tabular}[x]{@{}c@{}}0.015\\ $T_0=10$\\$Maxit=10$\end{tabular}\\ \hline

$\Cay(\mathbb{Z}_9,\{1,2,3,6,7,8\})$&9 &4&0.41&\begin{tabular}[x]{@{}c@{}}0.01\\ $n_p=10$\\$Maxit=50$\end{tabular}&\begin{tabular}[x]{@{}c@{}}0.017\\ $T_0=10$\\$Maxit=10$\end{tabular}\\ \hline

$\Cay(\mathbb{Z}_{20},\{1,3,17,19\})$&20 &7&0.85&\begin{tabular}[x]{@{}c@{}}0.54\\$n_p=20$\\$Maxit=250$\end{tabular}&\begin{tabular}[x]{@{}c@{}}0.286\\ $T_0=10$\\$Maxit=50$\end{tabular}\\ \hline
\end{tabular}
\end{center}
\caption{The general position number of some graphs and running times of the ILP, the GA, and the SA} 
\label{Tab1}
\end{table}

As can be seen from Table~\ref{Tab1}, SA was the fastest in six cases and GA in four cases. Both algorithms generally seem to be faster than ILP, which has never been the fastest.

\section*{Acknowledgments}

Sandi Klav\v zar were supported by the Slovenian Research and Innovation Agency (ARIS) under the grants P1-0297, N1-0355, and N1-0285.

\end{document}